\documentclass[12pt,letterpaper]{amsart}
\usepackage{euler, epic,eepic,latexsym, amssymb, amscd, amsfonts, xypic, floatflt, url, color, epsfig}

\newcommand{\tpoint}[1]{\vspace{3mm}\par \noindent \refstepcounter{subsection}{\bf \thesubsection.} 
  {\em #1. ---} }
\newcommand{\epoint}[1]{\vspace{3mm}\par \noindent \refstepcounter{subsection}{\thesubsection.} 
  {\em #1.} }
\newcommand{\bpoint}[1]{\vspace{3mm}\par \noindent \refstepcounter{subsection}{\bf \thesubsection.} 
  {\bf #1.} }

\newcommand{\bpf}{\noindent {\em Proof.  }}
\newcommand{\epf}{\qed \vspace{+10pt}}

\newcommand{\Z}{\mathbb{Z}}
\newcommand{\Zhat}{\hat{\mathbb{Z}}}

\newcommand{\kbar}{\overline{k}}

\newcommand{\G}{\mathbb{G}}

\newcommand{\proj}{\mathbb P}

\newcommand{\pmk}{\mathbb P^1_{k} - \{0,1,\infty \}}

\newcommand{\Gal}{\operatorname{Gal}}

\newcommand{\Path}{\operatorname{Path}}

\newcommand{\Hom}{\operatorname{Hom}}

\newcommand{\Spec}{\operatorname{Spec}}

\newcommand{\Aut}{\operatorname{Aut}}

\newcommand{\hidden}[1]{\footnote{Hidden:  #1}}
\renewcommand{\hidden}[1]{}

\title[Obstructions to $\pi_1$ sections]
{$n$-Nilpotent obstructions to $\pi_1$ sections of $\mathbb P^1 - \{0,1,\infty \}$ and Massey products}

\author[Wickelgren]{Kirsten Wickelgren}\thanks{The author is supported by an American Institute of Math Five Year Fellowship.}

\address{Dept. of Mathematics, Harvard}

\email{wickelgren@post.harvard.edu}

\date{June 11, 2011}



\begin{document}

\begin{abstract}

Let $\pi$ be a pro-$\ell$ completion of a free group, and let $G$ be a profinite group acting continuously on $\pi$. First suppose the action is given by a character. Then the boundary maps $\delta_n: H^1(G, \pi/[\pi]_n) \rightarrow H^2(G, [\pi]_n/[\pi]_{n+1})$ are Massey products. When the action is more general, we partially compute these boundary maps. Via obstructions of Jordan Ellenberg, this implies that $\pi_1$ sections of $\pmk$ satisfy the condition that associated $n^{th}$ order Massey products in Galois cohomology vanish. For the $\pi_1$ sections coming from rational points, these conditions imply that $\langle (1-x)^{-1}, x^{-1}, x^{-1}, \ldots, x^{-1} \rangle = 0$  where $x$ in $H^1(\Gal(\overline{k}/k), \Z_{\ell}(\chi))$ is the image of an element of $k^*$ under the Kummer map.

\end{abstract}

\maketitle

\section{Introduction} Grothendieck's section conjecture predicts that the rational points of a proper smooth hyperbolic curve $X$ over a number field $k$ are in natural bijection with the conjugacy classes of sections of the homotopy exact sequence for the \'etale fundamental group \begin{equation}\label{hes} \xymatrix{1 \ar[r] & \pi_1(X_{\overline{k}}) \ar[r] & \pi_1(X) \ar[r] & \pi_1(\Spec k)=\Gal(\overline{k}/k) \ar[r] & 1},\end{equation}  where the conjugacy class of a section $s: \Gal(\overline{k}/k) \rightarrow \pi_1(X)$ is the set of those sections $g \mapsto \gamma s(g) \gamma^{-1}$ where $\gamma$ is a fixed element of  $\pi_1(X_{\overline{k}})$. The phrase ``$\pi_1$ section" in the title refers to a section $s$ of $\pi_1(X) \rightarrow \pi_1(\Spec k)$. For a non-proper smooth hyperbolic curve, rational points ``at infinity" determine ``bouquets" of sections of (\ref{hes}) in bijection with $H^1(\Gal(\overline{k}/k), \Zhat(\chi))$ --see \cite[p. 2]{Popbipadicsc}. (Here, $\chi$ denotes the cyclotomic character, and $\Zhat(\chi^n)$ denotes $\Zhat$ with Galois action given by $\chi^n$.) More specifically, let $X$ be a smooth, geometrically integral curve over $k$ with negative Euler characteristic. Let $\overline{X}$ denote the smooth compactification of $X$. The section conjecture predicts that $$( \coprod_{(\overline{X}-X)(k)} H^1(\Gal(\overline{k}/k), \Zhat(\chi))) \coprod X(k)$$ is in bijection with the conjugacy classes of sections of (\ref{hes}) via a ``non-abelian Kummer map" discussed in \ref{kappadefn}.

 Consider the problem of counting the conjugacy classes of sections of (\ref{hes}). When (\ref{hes}) is split, this is equivalent to computing the pointed set $H^1(\Gal(\overline{k}/k), \pi_1(X_{\overline{k}}))$, which is difficult. In \cite{Ellenberg_2_nil_quot_pi}, Jordan Ellenberg suggested studying instead the image of $$H^1(\Gal(\overline{k}/k), \pi_1(X_{\overline{k}})) \rightarrow H^1(\Gal(\overline{k}/k),  \pi_1(X_{\overline{k}})^{ab})$$ by filtering $ \pi_1(X_{\overline{k}})$ by its lower central series. More specifically, let $\pi$ abbreviate $\pi_1(X_{\overline{k}})$, let $\pi^{ab}$ denote the abelianization of $\pi$, and let $[\pi]_n$ denote the $n^{th}$ subgroup of the lower central series (cf. \ref{notationsubsection2}). Ellenberg proposed successively computing the images of $$H^1(\Gal(\overline{k}/k), \pi/[\pi]_n) \rightarrow H^1(\Gal(\overline{k}/k), \pi^{ab})$$ via the boundary maps $$\delta_n: H^1(\Gal(\overline{k}/k), \pi/[\pi]_n) \rightarrow H^1(\Gal(\overline{k}/k), [\pi]_n/[\pi]_{n+1})$$ coming from the central extensions $$\xymatrix{1 \ar[r] & [\pi]_n/[\pi]_{n+1} \ar[r] & \pi/[\pi]_{n+1} \ar[r] & \pi/[\pi]_n \ar[r] &1} .$$ This paper makes two group cohomology computations relating $\delta_n$ to Massey products (Propositions \ref{nomonodromy_delta_n} and \ref{muJdeltanwithmonodromy}), and then applies them to study the $\pi_1$ sections of $\pmk$ (Corollary \ref{pmkpi1restcor}) and Massey products of elements of $H^1(\Gal(\overline{k}/k), \Z^{\Sigma}(\chi))$, where $\Sigma$ is the set of primes not dividing any integer less than the order of the Massey product, and $\Z^{\Sigma}$ denotes the pro-$\Sigma$ completion of $\Z$ (Corollary \ref{vanishMasseyCor}). 

More specifically, the content of this paper is as follows: section \ref{freegroupcharaction} computes $\delta_n : H^1(G, \pi/[\pi]_n) \rightarrow H^2(G, [\pi]_n/[\pi]_{n+1})$ when $\pi$ is a pro-$\Sigma$ completion of a free group with generators $\{ \gamma_1, \gamma_2, \ldots, \gamma_r \}$, where $\Sigma$ is any set of primes not dividing $n!$, and $G$ is a profinite group acting on $\pi$ by $$g \gamma_i = \gamma_i^{\chi(g)} $$ where $\chi: G \rightarrow \Z^{\Sigma}$ is a character. In this case, $\delta_n$ is determined by $n^r$ order $n$ Massey products -- see Proposition \ref{nomonodromy_delta_n}. The case of the trivial character with $G$ and $\pi$ replaced by discrete groups is essentially contained in \cite{Dwyer}. The generalization to non-trivial characters is not immediate; for instance, it depends on the existence of certain upper triangular matrices whose $N^{th}$ powers are given by multiplying the $i^{th}$ upper diagonal by $N^i$--see (\ref{hom_for_defining}) and Lemma \ref{Aldef}. (To obtain these matrices one must invert $n!$ or work with pro-$\Sigma$ groups. We do the later, although the former works as well.) This computation is then used to study $\delta_n$ where $\pi$ is as above for $r=2$, and $G$ is a group acting on $\pi$ by $$g(\gamma_1) =\gamma_1^{\chi(g)} $$ $$g(\gamma_2) =\mathfrak{f}(g)^{-1}\gamma_2^{\chi(g)}\mathfrak{f}(g) $$ where $\chi: G \rightarrow \Z^{\Sigma}$ is a character, and $\mathfrak{f}: G \rightarrow [\pi]_2$ is a cocycle taking values in the commutator subgroup of $\pi$. In this case, $\delta_n$ pushed forward by certain Magnus coefficients are Massey products. The Magnus coefficients in question are those associated to degree $n$ non-commuative monomials in two variables containing $n-1$ factors of one variable  -- see Proposition \ref{muJdeltanwithmonodromy}. This calculation imposes restrictions on the image of $$H^1(\Gal(\overline{k}/k), \pi/[\pi]_{n+1}) \rightarrow H^1(\Gal(\overline{k}/k), \pi^{ab})$$ for $X = \pmk$, $\pi = \pi_1(X_{\kbar})^{\Sigma}$. Identifying $H^1(\Gal(\overline{k}/k), \pi^{ab})$ with $H^1(\Gal(\overline{k}/k), \Z^{\Sigma}(\chi))^2$, these restrictions are that the image is contained in the subset of elements  $x_1 \times x_2$ such that the Massey products $\langle -x_{J(1)}, -x_{J(2)}, \ldots, -x_{J(n)} \rangle$ vanish for all $J: \{1,2,\ldots,n \} \rightarrow \{1,2\}$ which only assume the value $2$ once -- see Corollary \ref{pmkpi1restcor}. Corollary \ref{vanishMasseyCor} writes these restrictions for the $\pi_1$ sections coming from rational points and tangential points, and concludes that the $n^{th}$ order Massey products $$\langle x^{-1}, \ldots, x^{-1}, (1-x)^{-1}, x^{-1}, \ldots, x^{-1} \rangle \textrm { and } \langle x,\ldots x,-x, x, \ldots, x \rangle$$ vanish, where $x$ in $H^1(\Gal(\overline{k}/k), \Zhat(\chi))$ denotes the image of an element of $k^*$ under the Kummer map.  Much of this vanishing behavior was previously shown by Sharifi \cite{Sharifi}, who calculates Massey products of the form $\langle x,x,\ldots,x,y\rangle$ under certain hypotheses and using different methods --see remarks \ref{pi1restremarks} and \ref{Sharifivanishing}. Triple Massey products in Galois cohomology with restricted ramification are studied by Vogel in \cite{Vogel_thesis}. 

The first subsections of sections \ref{deltansection} and \ref{application_section} contain only well-known material. They are meant to be expository and to fix notation.

{\em Acknowledgments:} I wish to thank Romyar Sharifi for useful correspondence.

\section{$n^{th}$ order Massey products and $\delta_n$}\label{deltansection}

\bpoint{Notation}\label{notationsubsection2} 
For elements $g_1$, $g_2$ of $G$, let $[g_1,g_2] = g_1 g_2 g_1^{-1} g_2^{-1}$ denote the commutator.  For a profinite group $\pi$, let $\pi=[\pi]_1 \supset [\pi]_2 \supset [\pi]_3 \ldots$ denote the lower central series: $[\pi]_n$ is defined to be the closure of the subgroup generated by the elements of $[\pi, [\pi]_{n-1}]$.

\label{Z(chi)notation} For a (profinite) group $G$, a profinite abelian group $A$, and a (continuous) homomorphism $\chi: G \rightarrow \Aut(A)$, let $A(\chi)$ denote the associated profinite group with $G$ action. For example, if $A$ is a ring and $\chi$ is a homomorphism $G \rightarrow A^*$, then for any integer $n$, $A(\chi^n)$ is a profinite group with $G$ action.

 Let $\Sigma$ denote a set of primes (of $\Z$). For any group $G$, let $G^{\Sigma}$ denote the pro-$\Sigma$ completion of $G$, i.e. the inverse limit of all quotients of $G$ whose order divides a product of powers of primes in $\Sigma$.

\bpoint{Massey Products}\label{Masseydefn}

For a profinite group $G$ and a profinite abelian group $A$ with a continuous action of $G$, let $(C^*(G,A), D)$ be the complex of inhomogeneous cochains of $G$ with coefficients in $A$ as in \cite[I.2~p. 14]{coh_num_fields}.  For $c \in C^p(G, A)$ and $d \in C^q(G, A)$, let $c \cup d$ denote the cup product $c \cup d \in C^{p+q}(G, A \otimes A)$  $$(c \cup d)(g_1,\ldots,g_{p+q}) = c(g_1,\ldots,g_p) \otimes ((g_1\cdots g_p)d(g_{p+1},\ldots, g_{p+q})).$$ This product induces a well defined map on cohomology. If $A$ is a ring with $G$ action, then the action is given by a homomorphism $\chi: G \rightarrow A^*$. Recall the notation $A(\chi^n)$ defined in \ref{Z(chi)notation}. The $G$ equivariant multiplication map $A(\chi^n) \otimes A(\chi^m) \rightarrow A(\chi^{n+m})$ induces cup products $$C^p(G,A(\chi^n)) \otimes C^q(G,A(\chi^m)) \rightarrow C^{p+q}(G, A(\chi^{n+m})) $$ $$H^p(G,A(\chi^n)) \otimes H^q(G,A(\chi^m)) \rightarrow H^{p+q}(G, A(\chi^{n+m})) .$$ 

For a profinite group $Q$, no longer assumed to be abelian, the set of continuous functions $G \rightarrow Q$ is denoted $C^1(G,Q)$. An element $s$ of $C^1(G,Q)$ such that $s(g_1 g_2) = s(g_1) g_1 s(g_2)$ is a {\em cocycle} or {\em twisted homomorphism}. $H^1(G,Q)$ is defined as equivalence classes of cocycles in the usual manner (cf. \cite[VII Appendix]{serre:localfields}). 

\epoint{Definition}\label{Massey_prod_def} Let $t_1,\ldots, t_n$ be elements of $H^1(G,A(\chi)).$ The $n^{th}$ order Massey product of the ordered $n$-tuple $(t_1,\ldots, t_n)$ is defined if there exist $ T_{ij}$ in $C^1(G,A(\chi^{j-i}))$ for $i,j$ in $\{1,2,\ldots, n+1 \}$ such that $i<j$ and $(i,j) \neq (1,n+1)$ satisfying

\begin{itemize}
\item $T_{i,i+1}$ represents $t_i$.
\item $D T_{ij} = \sum_{p=i+1}^{j-1} T_{ip} \cup T_{pj}$ for $i+1<j$
\end{itemize}

$T$ is called a {\em defining system}. The {\em Massey product relative to $T$} is defined by $$\langle t_1,\ldots t_{n-1} \rangle_T = \sum_{p=2}^{n} T_{1p} \cup T_{p,n+1} .$$ 

\epoint{Massey products and unipotent matrices} Let $U_{n+1}$ denote the multiplicative group of $(n+1) \times (n+1)$ upper triangular matrices with coefficients in $A$ whose diagonal entries are $1$. (``U" stands for unipotent-- not unitary.) Let $a_{ij}$ be the function taking a matrix to its $(i,j)$-entry. $U_{n+1}$ inherits an action of $G$ by $a_{ij} (g M )= \chi(g)^{j-i} a_{ij}M$. We have a $G$ equivariant inclusion $A(\chi^n) \rightarrow U_{n+1}$ sending $a$ in $A$ to the matrix with $a$ in the $(1,n)$-entry, and with all other off diagonal matrix entries $0$.  This inclusion gives rise to a central extension

\begin{equation}
\label{A_U_barU}
1 \rightarrow A(\chi^n) \rightarrow U_{n+1} \rtimes G \rightarrow \overline{U}_{n+1} \rtimes G \rightarrow 1.
\end{equation}

\noindent where $\overline{U}_{n+1}$ is defined as the quotient $U_{n+1}/A(\chi^n)$. 

The element of $H^2(\overline{U}_{n+1} \rtimes G, A(\chi^n))$ classifying (\ref{A_U_barU}) is an order $n$ Massey product. (\label{H2kelement}See \cite[IV \S 3]{Brown_coh_groups} for the definition of the element of $H^2$ classifying a short exact sequence of groups; to apply the same discussion to profinite groups, one needs continuous sections of profinite quotient maps. For this, see \cite[Prop 2.2.2]{Profinite_Groups}.) Note that $-a_{i,j}$ determines an element of $C^1(\overline{U}_{n+1} \rtimes G, A(\chi^{j-i}))$. As $(i,j)$ ranges through the set of pairs of elements of $\{1,2,\ldots, n+1 \}$ such that $i<j$ and $(i,j) \neq (1,n+1)$, $-a_{i,j}$ is a defining system for $(-a_{1,2}, -a_{2,3}, \ldots, -a_{n,n+1})$. The element of $H^2(\overline{U}_{n+1} \rtimes G, A(\chi^n))$ classifying (\ref{A_U_barU}) is $\langle -a_{1,2}, -a_{2,3}, \ldots, -a_{n,n+1} \rangle$, where the Massey product is taken with respect to the defining system $-a_{i,j}$. This follows immediately from the definition of matrix multiplication. 

\bpoint{Magnus embedding}\label{Magnus_embedding_recall} For later use, we recall some well known properties of the Magnus embedding. Let $F$ denote the free group on the $r$ generators $\gamma_i$, $i = 1,\ldots, r$. For any ring $A$, let $A \langle \langle z_1,\ldots, z_r \rangle \rangle$ be the ring of associative power series in the non-commuting variables $z_1, \ldots,z_r$ with coefficients in $A$. Let $A \langle \langle z_1,\ldots,z_r \rangle \rangle^{(1,\times)}$ denote the subgroup of the multiplicative group of units of $A \langle \langle z_1,\ldots,z_r \rangle \rangle$ consisting of power series with constant coefficient $1$. The {\em Magnus embedding} is defined $$F \rightarrow \Z \langle \langle z_1,\ldots,z_r \rangle \rangle^{(1,\times)}$$ by $x_j \mapsto 1+ z_j$ for all $j$. 

Since $\Z^{\Sigma} \langle \langle z_1,\ldots,z_r \rangle \rangle^{(1,\times)}$ is pro-$\Sigma$, $F \rightarrow \Z \langle \langle z_1,\ldots,z_r \rangle \rangle^{(1,\times)}$ gives rise to a commutative diagram $$\xymatrix{ F^{\Sigma} \ar[r] & \Z^{\Sigma} \langle \langle z_1,\ldots,z_r \rangle \rangle^{(1,\times)}\\
F \ar[u] \ar[r] & \ar[u] \Z \langle \langle z_1,\ldots,z_r \rangle \rangle^{(1,\times)}}.$$ Let $J:\{1,\ldots,n\} \rightarrow \{1,\ldots,r\}$ be any function. The degree $n$ monomial $ z_{J(1)} \cdots z_{J(n)}$ determines the {\em Magnus coefficient} $\mu_J: F^{\Sigma} \rightarrow \Z^{\Sigma}$ (or $\mu_J : F \rightarrow \Z$ ) given by taking an element of $F^{\Sigma}$ to the coefficient of $ z_{J(1)} \cdots z_{J(n)}$ in its image under the Magnus embedding. It is well known that $\mu_J (\gamma) =0$ for $\gamma \in [F]_m$ and $m>n \geq 1$ (see \cite[\S 5.5, ~ Cor.~5.7]{Magnus_Karrass_Solitar}), and it follows by continuity that $\mu_J (\gamma) =0$ for $\gamma$ in $F^{\Sigma}$ and $m>n \geq 1$. 

The {\em Lie elements} of $\Z \langle \langle z_1,\ldots,z_r \rangle \rangle$ are the elements in the image of the Lie algebra map $\zeta_i \mapsto z_i$ from the free Lie algebra over $\Z$ on $r$ generators $\zeta_i$ to $\Z \langle \langle z_1,\ldots,z_r \rangle \rangle$, where $\Z \langle \langle z_1,\ldots,z_r \rangle \rangle$ is considered as a Lie algebra with bracket $[z,z'] = z z' - z' z$.  It is well known that the Magnus embedding induces an isomorphism from $[F]_n/[F]_{n+1}$ to the homogeneous degree $n$ Lie elements of $\Z \langle \langle z_1,\ldots,z_r \rangle \rangle$ \cite[\S 5.7,~Cor.~5.12(i)]{Magnus_Karrass_Solitar}. The Lie basis theorem \cite[\S 5.6,~Thm.~5.8(ii)]{Magnus_Karrass_Solitar} implies that the inclusion of the Lie elements of degree $n$ into all the degree $n$ elements of $\Z \langle \langle z_1,\ldots,z_r \rangle \rangle$ is a direct summand. It follows that $$\xymatrix{
[F]_n/[F]_{n+1} \ar[rr]^{\oplus_J \mu_J}&& \oplus_J \Z } $$ is the inclusion of a (free) direct summand. By definition of $\mu_J$, we have the commutative diagram \begin{equation}\label{UJcomdiag}\xymatrix{ [F^{\Sigma}]_n/[F^{\Sigma}]_{n+1} \ar[rr]^{\oplus_J \mu_J} && \oplus_J \Z^{\Sigma} \\
[F]_n/[F]_{n+1} \ar[u] \ar[rr]^{\oplus_J \mu_J}&& \oplus_J \Z \ar[u]} \end{equation} where the direct sums are taken over all functions $$J:\{1,\ldots,n\} \rightarrow \{1,\ldots,r\}.$$ 

We claim that the top horizontal morphism in (\ref{UJcomdiag}) is the inclusion of a direct summand of the form $\oplus \Z^{\Sigma}$, and that the left vertical morphism is the pro-$\Sigma$ completion. To see this: note that since $[F]_n/[F]_{n+1}$ is a free $\Z$ submodule of $\oplus_J \Z$, we have a commutative diagram $$\xymatrix{ ([F]_n/[F]_{n+1})^{\Sigma} \ar[rr] && \oplus_J \Z^{\Sigma} \\
[F]_n/[F]_{n+1} \ar[u] \ar[rr]^{\oplus_J \mu_J}&& \oplus_J \Z \ar[u]} $$ where the bottom horizontal morphism is the inclusion of a direct summand which is a free $\Z$ module, the top horizontal morphism is the inclusion of a direct summand which is a free $\Z^{\Sigma}$ module, and both vertical maps are pro-$\Sigma$ completions. The map  $([F]_n/[F]_{n+1})^{\Sigma} \rightarrow \oplus_J \Z^{\Sigma}$ factors through $[F^{\Sigma}]_n/[F^{\Sigma}]_{n+1}$ by the universal property of pro-$\Sigma$ completion, and it follows that $([F]_n/[F]_{n+1})^{\Sigma}\rightarrow [F^{\Sigma}]_n/[F^{\Sigma}]_{n+1} $ is injective. Since $[F]_n$ has dense image in  $[F^{\Sigma}]_n$, and since the image of a compact set under a continuous map to a Hausdorff topological space is closed, we have that $([F]_n/[F]_{n+1})^{\Sigma} \rightarrow [F^{\Sigma}]_n/[F^{\Sigma}]_{n+1} $ is surjective. Thus $([F]_n/[F]_{n+1})^{\Sigma} \rightarrow [F^{\Sigma}]_n/[F^{\Sigma}]_{n+1} $ is an isomorphism of profinite groups (because a continuous bijection between compact Hausdorff topological spaces is a homeomorphism). From this it also follows that the top horizontal morphism in (\ref{UJcomdiag}) is the inclusion of a direct summand of the form $\oplus \Z^{\Sigma}$.

\bpoint{$\delta_n$ for the free pro-$\Sigma$ group with action via a character}\label{freegroupcharaction} Let $n$ be a positive integer and let $\Sigma$ be the set of primes not dividing $n!$. Let $\pi$ be the pro-$\Sigma$ completion of the free group on the generators $\{ \gamma_1, \gamma_2, \ldots, \gamma_r \}$. Let $G$ be any profinite group and let $\chi: G \rightarrow (\Z^{\Sigma})^*$ be a (continuous) character of $G$. Let $G$ act on $\pi$ via $g \gamma_i = \gamma_i^{\chi(g)}$. Then the map $$\delta_n: H^1(G, \pi/[\pi]_{n}) \rightarrow H^2(G, [\pi]_n/[\pi]_{n+1})$$ is given by $n^{th}$ order Massey products in the following manner:

Recall that $U_{n+1}$ denotes the group of $(n+1) \times (n+1)$ upper triangular matrices with diagonal entries equal to $1$, that $a_{i,j}: U_{n+1} \rightarrow \Z^{\Sigma}$ denotes the $(i,j)^{th}$ matrix entry, and that $U_{n+1}$ inherits a $G$-action defined by $a_{i,j} (g ( M )) = \chi(g)^{j-i} a_{i,j} (M)$ for all $M$ in $U_{n+1}$ (see \ref{Masseydefn}). 

For each $J: \{1,2,\ldots, n \} \rightarrow \{1,2,\ldots,r\}$, let $\varphi_J: \pi \rightarrow U_{n+1}$ be the homomorphism defined \begin{equation}\label{hom_for_defining} a_{i, j} \varphi_J (\gamma_k) = 
 \begin{cases}
 \frac{1}{l!} & j = i+ l, l>0 \mbox{ and } k = J(v)\mbox{ for all }i \leq v < i+l \\
 1 & $j=i$ \\
 0 & \mbox{otherwise}
 \end{cases} \end{equation}
 
It is a straightforward consequence of the following lemma that $\varphi_J$ is $G$ equivariant:

 \tpoint{Lemma}\label{Aldef}{\em Let $A_l$ be the matrix in $U_{l+1}$ defined by $a_{i, i+j}A_l = \frac{1}{j!}$ for $j>0$. Then for all positive integers $N$, $a_{i,i+j}(A_l^N )= N^j a_{i,i+j}(A_l)$.}
 
\bpf
By induction on $l$. For $l=1$, the lemma is clear. By induction and symmetry, it is sufficient to check that $a_{1,l+1}(A_l^N )= N^l a_{1,l+1}(A_l)$. Now induct on $N$, so in particular, $a_{1,1+j}(A_l^{N-1}) = (N-1)^j \frac{1}{j!}$ for $j=0,\ldots,l$. Thus $$a_{1,l+1}(A_l^N )= \sum_{j=0}^{l} a_{1, 1+j}(A_l^{N-1}) \frac{1}{(l-j)!} = $$ $$\sum_{j=0}^{l} (N-1)^j \frac{1}{j!} \frac{1}{(l-j)!} = ((N-1) + 1)^l \frac{1}{l!},$$ completing the proof.
\epf

Since $[U_{n+1}]_{n+1} = 1$ and $[\overline{U}_{n+1}]_n = 1$, $\varphi_J$ descends to $G$-equivariant homomorphisms $$\pi/[\pi]_{n+1} \rightarrow U_{n+1},$$  $$\pi/[\pi]_n \rightarrow \overline{U}_{n+1},$$ $$[\pi]_n/[\pi]_{n+1} \rightarrow \Z^{\Sigma}$$ which we also denote by $\varphi_J$.

The basis $\{\gamma_1, \gamma_2, \ldots, \gamma_r \}$ determines an isomorphism $$\pi/[\pi]_2 \cong \Z^{\Sigma}(\chi)^r,$$ and therefore an isomorphism $H^1(G, \pi^{ab}) \cong H^1(G, \Z^{\Sigma}(\chi))^r$. An element $x$ of $H^1(G, \pi/[\pi]_{n})$ projects to an element of $H^1(G, \pi^{ab})$. Let $x_1 \oplus \ldots \oplus x_r$ in $H^1(G, \Z^{\Sigma}(\chi))^r$ denote the image of the projection. 

Note that applying $a_{i,i+1} \varphi_J$ to a cocycle $x: G \rightarrow \pi/[\pi]_{n}$ produces a cocycle representing $x_{J(i)}$. Furthermore, $$\{ - a_{i,j} \varphi_J x : i< j, (i,j) \neq (1,n+1) \}$$ is a defining system for the Massey product $\langle -x_{J(1)}, -x_{J(2)}, \ldots, -x_{J(n)} \rangle$.

\tpoint{Proposition}\label{nomonodromy_delta_n}{\em For any cocycle $x: G \rightarrow \pi/[\pi]_n$, let $[x]$ denote the corresponding element of $H^1(G, \pi/[\pi]_n)$. Then $\delta_n ([x]) = 0$ if and only if $\langle - x_{J(1)}, - x_{J(2)}, \ldots,- x_{J(n)} \rangle = 0$ for every $J: \{1,2,\ldots, n \} \rightarrow \{1,2,\ldots,r\}$, where the Massey product is taken with respect to the defining system $\{ - a_{i,j} \varphi_J x : i< j, (i,j) \neq (1,n+1) \}$.}

\bpf
Choose $J: \{1,2,\ldots, n \} \rightarrow \{1,2,\ldots,r\}$. $\varphi_J$ induces a commutative diagram \begin{equation}\label{Utopi}
\xymatrix{1 \ar[r] & \Z^{\Sigma}(\chi^n) \ar[r] & U_{n+1} \rtimes G \ar[r] & \overline{U}_{n+1} \rtimes G \ar[r] & 1 \\
1 \ar[r] & [\pi]_n/[\pi]_{n+1} \ar[r] \ar[u] &\pi/[\pi]_{n+1} \rtimes G \ar[r] \ar[u] & \pi/[\pi]_{n} \rtimes G \ar[u] \ar[r] & 1}
 \end{equation} 
All the vertical morphisms in (\ref{Utopi}) will be denoted by $\varphi_J$. Let $\kappa$ denote the element of $H^2(\pi/[\pi]_{n} \rtimes G,  [\pi]_n/[\pi]_{n+1})$ classifying the bottom horizontal row, and let $\kappa'$ denote the element of $H^2(\overline{U_{n+1}} \rtimes G,\Z^{\Sigma}(\chi^n))$  classifying the top horizontal row (c.f. \ref{H2kelement}). The morphism of short exact sequences (\ref{Utopi}) gives the equality $(\varphi_J)_* \kappa = \varphi_J^* \kappa'$ in $H^2(\pi/[\pi]_{n} \rtimes G,\Z^{\Sigma}(\chi^n))$. 

Choose a cocycle $x: G \rightarrow \pi/[\pi]_n$. Let $x \rtimes \operatorname{id}: G \rightarrow \pi/[\pi]_{n} \rtimes G$ denote the homomorphism $g \mapsto x(g) \rtimes g$ induced by the twisted homomorphism $x$. Then, $\delta_n([x]) = (x \rtimes \operatorname{id})^* \kappa$. Since $(\varphi_J)_* \kappa = \varphi_J^* \kappa'$, we have that $(\varphi_J)_* \delta_n([x]) = (\varphi_J \circ (x \rtimes \operatorname{id}))^*  \kappa'$. By \ref{Masseydefn},  $(\varphi_J \circ (x \rtimes \operatorname{id}))^*  \kappa'$ is the Massey product $\langle - x_{J(1)},- x_{J(2)}, \ldots,- x_{J(n)} \rangle$ computed with the defining system $\{- a_{i,j} \varphi_J x : i< j, (i,j) \neq (1,n+1) \}$. 

It is therefore sufficient to see that $$\oplus_J (\varphi_J)_*: H^2(G, [\pi]_n/[\pi]_{n+1}) \rightarrow H^2(G, \oplus_J \Z^{\Sigma}(\chi^n))$$ is injective. This follows from a result of Dwyer: let $\mu_J$ denote the Magnus coefficient as in \ref{Magnus_embedding_recall}. By \cite[Lem 4.2]{Dwyer}, $\mu_J (\gamma)= \varphi_J (\gamma)$ for any element $\gamma$ in the free group generated by the $\gamma_i$, and the equality $\mu_J = \varphi_J$ for any element of $\pi$ follows by continuity. Thus the map $\oplus_J (\varphi_J): [\pi]_n/[\pi]_{n+1} \rightarrow \oplus_J \Z^{\Sigma}(\chi^n)$ is the split injection $\oplus_J \mu_J$ induced by the homogeneous degree $n$ piece of the Magnus embedding -- see (\ref{UJcomdiag}) in \ref{Magnus_embedding_recall}.
 \epf

Thus, if the element $x_1 \oplus x_2 \oplus \ldots \oplus x_r$ of $H^1(G, \pi^{ab})$ lifts to $x$ in $$H^1(G, \pi/[\pi]_{n+1}),$$ all the order $n$ Massey products $$\langle -x_{J(1)}, - x_{J(2)}, \ldots, - x_{J(n)}\rangle = (\varphi_J)_* \delta_n(x)$$ vanish. Furthermore, if the vanishing of the order $n$ Massey products occurs with respect to defining systems which are compatible in the sense of Proposition \ref{nomonodromy_delta_n}, the converse holds as well.

\bpoint{Partial computation of $\delta_n$ for the free pro-$\Sigma$ group with action determined by a character and conjugation by a cocyle valued in the commutator subgroup}\label{monodromy_deltan_bpoint} Choose a positive integer $n$, and let $\Sigma$ denote the set of primes not dividing $n!$. Let $\pi$ be the pro-$\Sigma$ completion of the free group on two generators $\{ \gamma_1, \gamma_2\}$. Let $\chi: G \rightarrow (\Z^{\Sigma})^*$ be a (continuous) character of a profinite group $G$. Let $G$ act on $\pi$ via \begin{align}
\label{G_actionpiP1-} 
 g(\gamma_1) &=\gamma_1^{\chi(g)} \\ g(\gamma_2) &=\mathfrak{f}(g)^{-1}\gamma_2^{\chi(g)}\mathfrak{f}(g), \nonumber
 \end{align}
 
\noindent where $\mathfrak{f}: G \rightarrow [\pi]_2$ is a cocyle. For instance, the Galois action on the pro-$\Sigma$ \'etale fundamental group of $\pmk$ has this form with respect to an appropriate base point. (See, for instance, \cite{Ihara_GT}. This situation will be considered in section \ref{application_section}.)

Then there are obstructions to $\delta_n = 0$ given by order $n$ Massey products:

Choose $i_0$ in $ \{1,2,\ldots,n\}$ and let $J: \{1,2,\ldots,n\} \rightarrow \{1,2\}$ be the function $J(i_0) = 2$, $J(j) = 1$ for $j\neq i_0$. Let $\varphi_J: \pi \rightarrow U_{n+1}$ be the homomorphism given by equation (\ref{hom_for_defining}) in \ref{freegroupcharaction}. The next two lemmas are used to show that $\varphi_J$ is $G$-equivariant. 

\tpoint{Lemma}\label{Ui0j0normalcommutative}{\em Let $$U_{i_0, j_0} = \{ M \in U_{n+1}: a_{ij}(M) = 0 \textrm{ for } i \neq j \textrm{ unless }  i \leq i_0 \textrm{ and } j \geq j_0 \}$$ Then $U_{i_0, j_0}$ is a normal subgroup of $U$ which is commutative for $i_0 < j_0$.}

\bpf
It is straightforward to see that $U_{i_0, j_0}$ is a subgroup. (Indeed, it suffices to note that for $M_1$,$M_2$ in $U_{i_0, j_0}$, we have $a_{ij} ((M_1 -1)(M_2 -1)) = 0$ for $i>i_0$ or $j<j_0$; for instance, this implies that $U_{i_0, j_0}$ is closed under inverses because $M^{-1} = 1 + \sum_{k \geq 1} (-1)^k (M-1)^k$.) 

To see that $U_{i_0, j_0}$ is normal, choose $Z$ in $U_{n+1}$ and $M$ in $U_{i_0, j_0}$. Note that $$a_{ij}(Z (M-1) Z^{-1}) = \sum_{k,k'} a_{ik}(Z_{ik}) a_{k k'} (M-1) a_{k'j} (Z^{-1})$$ which is only non-zero if there exists $k$ and $k$ such that $ i \leq k \leq k' \leq j$, $k \leq i_0$, and $k' \geq j_0$. This can only occur for  $i \leq i_0$ and $j \geq j_0$. So, $U_{i_0, j_0}$ is normal.

Suppose that $i_0 < j_0$, and let $M_1$, $M_2$ be in $U_{i_0, j_0}$. To see that $U_{i_0, j_0}$ is commutative, it suffices to see that $(M_1 -1) (M_2 -1) = 0$. To see this equality, note that for all $i,j,k$, we have $a_{i k}(M_1 - 1) a_{k j} (M_2 -1) = 0$, because if $k< j_0$, then $a_{i k}(M_1 - 1) =0$, and if $k \geq j_0$, then $k > i_0$, whence $ a_{k j} (M_2 -1) = 0$.
\epf

\tpoint{Lemma}{\em $\varphi_J(\gamma_2)$ commutes with any element of $\varphi_J([\pi]_2)$.}

\bpf
Let $X = \varphi_J(\gamma_1)$, $Y = \varphi_J(\gamma_2)$, and $\varpi$ be the closure of the subgroup generated by $X$ and $Y$. By Lemma \ref{Ui0j0normalcommutative}, it is sufficient to show that $[\varpi]_2$ is contained in $U_{i_0, i_0 + 1}$. 

$[\varpi]_2$ is topologically generated by elements of the form $$[\cdots [[X,Y], Z_1], Z_2, \ldots], Z_k]$$ where $Z_i$ is either $X$ or $Y$ and $k = 0,1,\ldots$. By Lemma  \ref{Ui0j0normalcommutative}, if $W$ is in $U_{i_0, i_0 + 1}$, so is $[W,Z]$ for any $Z$ in $\varpi$. Since $Y$ is in $U_{i_0, i_0 + 1}$, it follows that $[\cdots [[X,Y], Z_1], Z_2, \ldots], Z_k]$ is as well.
\epf

In particular, $\varphi_J(g(\gamma_2)) = \varphi_J(\gamma_2)^{\chi(g)}$, so $\varphi_J$ is $G$-equivariant by \ref{freegroupcharaction}. 

Since $\varphi_J$ is $G$-equivariant, we have the commutative diagram (\ref{Utopi}). Choose a cocycle $x: G \rightarrow \pi/[\pi]_n$. $(\varphi_J)_* \delta_n([x]) $ is the Massey product $\langle - a_{1,2} \varphi_J x ,- a_{2,3} \varphi_J x , \ldots,- a_{n,n+1} \varphi_J x \rangle$ computed with the defining system $\{- a_{i,j} \varphi_J x : i< j, (i,j) \neq (1,n+1) \}$ by \ref{Masseydefn} (see the proof of Proposition \ref{nomonodromy_delta_n}).

Note that $\{\gamma_1,\gamma_2\}$ is a $\Z^{\Sigma}(\chi)$ basis for $\pi/[\pi]_2$, giving an isomorphism $H^1(G, \pi^{ab}) \cong H^1(G, \Z^{\Sigma}(\chi))^2$. As above, an element $x$ of $H^1(G, \pi/[\pi]_{n})$ projects to an element of $H^1(G, \pi^{ab})$. Let $x_1 \oplus x_2$ in $H^1(G, \Z^{\Sigma}(\chi))^2$ denote the image of the projection. Note that $- a_{j,j+1} \varphi_J x = x_{J(j)}$. We have therefore shown:

\tpoint{Proposition}\label{muJdeltanwithmonodromy}{\em Let $x: G \rightarrow \pi/[\pi]_n$ be a cocycle, and let $[x]$ denote the corresponding cohomology class. If $\delta_n([x]) = 0$, then $\langle - x_{J(1)}, - x_{J(2)}, \ldots, - x_{J(n)} \rangle = 0$ where this Massey product is taken with respect to the defining system $\{ - a_{i,j} \varphi_J x : i< j, (i,j) \neq (1,n+1) \}$ defined above.}

\epoint{Remark} As in \ref{freegroupcharaction}, the Massey product in Proposition \ref{muJdeltanwithmonodromy} is $\mu_J \delta_n$, where $\mu_J$ is the Magnus coefficient defined in \ref{Magnus_embedding_recall}. In other words, Proposition \ref{muJdeltanwithmonodromy} computes $\mu_J \delta_n$ for all functions $J$ which only assume the value $2$ once.

\section{Application to $\pi_1$ sections of punctured $\proj^1$ and Massey products in Galois cohomology}\label{application_section}

\bpoint{Notation}\label{Notation3} Let $k$ be a field of characteristic $0$ and let $\overline{k}$ be an algebraic closure of $k$. Let $G_k= \Gal(\overline{k}/k)$ denote the absolute Galois group of $k$.

\bpoint{Review of the \'etale fundamental group}\label{etpi1rev} A geometric point $b$ of a scheme $X$ (i.e. a map $b: \Spec \Omega \rightarrow X$ where $\Omega$ is an algebraically closed field) determines a functor from the finite \'etale covers of $X$ to the category of sets, called a {\em fiber functor}. Given two geometric points $b_1$, $b_2$, define $\Path(b_1,b_2)$ to be the natural transformations from the fiber functor associated to $b_1$ to the fiber functor associated to $b_2$. $\Path(b_1,b_2)$ naturally has the structure of a profinite set. Path composition will be in ``functional order," so given $\wp_1$ in $\Path(b_1,b_2)$ and $\wp_2$ in $\Path(b_2,b_3)$, we have $\wp_2 \wp_1$ in $\Path(b_1,b_3)$. The \'etale fundamental group $\pi_1(X,b)$ is the profinite group $\Path (b,b)$ (see \cite{sga1} \cite{Mezard}).

Suppose that $X$ is defined over a field $k$. Let $\overline{k}$ denote a fixed algebraic closure of $k$, and let $X_{\overline{k}} = X \times_{\Spec k} \Spec \overline{k}$ denote the base change of $X$ to $\overline{k}$. A rational point $\Spec k \rightarrow X$ gives rise to a geometric point $\Spec \overline{k} \rightarrow X_{\overline{k}}$ of $X_{\overline{k}}$, and there is a natural action of $G_k$ on the associated fiber functor, induced by the commutative diagram $$\xymatrix{ \Spec \overline{k} \ar[r]^{g} \ar[d] & \Spec \overline{k} \ar[d] \\ X_{\overline{k}} \ar[r]^{g} & X_{\overline{k}} } $$ where $g$ is any element of $G_k$. 

Now suppose that $X$ is a smooth, geometrically connected curve over $k$. Let $\overline{X}$ denote the smooth compactification, and let $x: \Spec k \rightarrow \overline{X}$ be a rational point. The completed local ring of $\overline{X}$ at the image of $x$ is isomorphic to $k[[\epsilon]]$ and the choice of such an isomorphism gives a map $\Spec k((\epsilon)) \rightarrow X$, where $k((\epsilon))$ denotes the field of Laurent power series. Such a map will be called a {\em rational tangential point}. To a rational tangential point, we can naturally associate a map $\Spec \overline{k}((\epsilon)) \rightarrow X_{\overline{k}}$. Since $k$ is characteristic $0$, the field of Puiseux series $\cup_{n \in \Z_{> 0}} \overline{k} ((\epsilon^{1/n}))$ is algebraically closed. Embedding  $\overline{k}((\epsilon))$ in  $\cup_{n \in \Z_{> 0}} \overline{k} ((\epsilon^{1/n}))$ in the obvious manner allows us to associate to a rational tangential point a corresponding geometric point $\Spec \cup_{n \in \Z_{> 0}} \overline{k} ((\epsilon^{1/n})) \rightarrow X_{\overline{k}}$ and fiber functor. There is a $G_k$ action on this fiber functor given by the previous commutative diagram with the field of Puiseux series replacing $\overline{k}$, and where $g$ in $G_k$ is taken to act on the field of Puiseux series via the action on the $\overline{k}$ coefficients. Tangential points are discussed in greater generality and more detail in \cite[\S 15]{Deligne} and \cite{Nakamura}.

\epoint{Example}\label{tgtvecA1} Let $U$ be an open subset of $\mathbb{A}^1_{k} = \Spec k[z]$. A tangent vector of $\mathbb{A}^1_{k}$ $$\Spec k[\epsilon]/\langle \epsilon^2 \rangle \rightarrow \mathbb{A}^1_{k}$$ $$z \mapsto a + v \epsilon $$ where $a$ in $k$, $v$ in $k^*$, gives a rational tangential point $\Spec k((\epsilon)) \rightarrow U$ of $U$ by $z \mapsto a + v \epsilon $.

By a {\em rational base point}, we will mean either a rational point or a rational tangential point, and by a slight abuse of notation, {\em rational base point} will also refer to the geometric points given above and their fiber functors. 

\epoint{Construction: ``non-abelian Kummer map"}\label{kappadefn} Let $X$ be a smooth curve over a field $k$. Let $X^{bp}(k)$ denote the set of rational base points of $X$, and assume that $X^{bp}(k) \neq \emptyset$. Fix $b$ in $X^{bp}(k)$. There is a ``non-abelian Kummer map" based at $b$ $$\kappa_b: X^{bp}(k) \rightarrow H^1(G_k, \pi_1(X_{\overline{k}}, b))$$ defined as follows: for $x$ in $X^{bp}(k)$, choose $\wp$ in $\Path(b, x)$ and define a $1$-cocycle $\kappa_{(b, \wp)} (x): G_k \rightarrow \pi_1(X_{\overline{k}}, b)$ by \begin{equation}\label{eqkappadef}\kappa_{(b, \wp)} (x) (g) = \wp^{-1} (g \wp).\end{equation}  \label{kappagammadefn} The cohomology class of this cocycle is independent of the choice of $\wp$ and $\kappa_b (x)$ is defined to be this cohomology class. When the base point is clear, $\kappa_b$ will also be denoted by $\kappa$.

Note that associated to a rational tangential point of $X$, there is a tangent vector $$\Spec k[\epsilon]/\langle \epsilon^2 \rangle \rightarrow \overline{X}$$ of the smooth compactification (see the above definition of a rational tangential point; the tangent vector is induced by the chosen isomorphism of $k[[\epsilon]]$ with the completed local ring of $\overline{X}$). It is not difficult to check that the images under $\kappa_b$ of two rational tangential points with the same tangent vector are equal (see \cite[p 6]{PIA}).

\epoint{Kummer map in Galois cohomology} Let $k$ be a field of characteristic $0$, and choose an isomorphism of the roots of unity in $\overline{k}$ with $\Zhat(\chi)$, where $\chi$ denotes the cyclotomic character. The short exact sequence $$\xymatrix{1 \ar[r] & \Z/m(\chi) \ar[r] & \overline{k}^* \ar[r]^{x \mapsto x^m} & \overline{k}^* \ar[r] & 1}$$ of $G_k$ modules gives a boundary map $k^* \rightarrow H^1(G_k, \Z/m(\chi))$. Letting $m$ vary gives the Kummer map $$k^* \rightarrow H^1(G_k, \Z^{\Sigma}(\chi)) $$ We will adopt the notational convention that an element of $k^*$ will also denote the corresponding class in $H^1(G_k, \Z^{\Sigma}(\chi))$. 

\tpoint{Lemma}\label{tgtl_pts_0_Gm}{\em For $\G_m$ based at the rational point $1$, $\kappa(x) = x$ and $\kappa(0 + v \epsilon) = v$ for all $x,v$ in $k^*$.} 

We sketch of a proof of this well-known fact.

\bpf
 The connected finite \'etale covers of $\G_{m, \overline{k}} = \Spec \overline{k}[z, \frac{1}{z}]$ are $\xymatrix{\G_{m, \overline{k}} \ar[rr]^{z \mapsto z^n} &&\G_{m, \overline{k}}}$ for $n$ in $\Z_{>0}$. Let $\mathcal{F} (0+ v \epsilon, n)$ denote the fiber of $z \mapsto z^n$ over (the geometric point associated to) $0 + v \epsilon$, where $0 + v \epsilon$ denotes the tangent vector $ k[\epsilon]/\langle \epsilon^2 \rangle \rightarrow \Spec k[z] $ given by $z \mapsto v \epsilon$ (cf. \ref{tgtvecA1}). Note that the $n^{th}$ roots of $v$ are in bijection with $\mathcal{F} (0+ v \epsilon, n)$; namely, an $n^{th}$ root $\sqrt[n]{v}$ of $v$ gives a map $\overline{k}[z^{1/n}, \frac{1}{z}] \rightarrow \cup_{n' \in \Z_{>0}} \overline{k} (( \epsilon^{1/n'}))$ which is tautologically a point of this fiber. Define $\mathcal{F} (1, n)$ similarly, and note that there is an identification of $\mathcal{F} (1, n)$ with the $n^{th}$ roots of unity in $\overline{k}$. A choice $\{ \sqrt[n]{v} : n \in \Z\}$ of compatible $n^{th}$ roots of unity of $v$ gives rise to $\wp$ in $\Path (1, 0 + v \epsilon)$; $\wp$ is the natural transformation such that the induced map $\mathcal{F} (1, n) \rightarrow \mathcal{F} (0+ v \epsilon, n)$ takes  $1$ to $\sqrt[n]{v}$. It follows that $g \wp$ takes $g 1$ to $g \sqrt[n]{v}$, from which we see that $\kappa(0 + v \epsilon) = v$. The equality $\kappa(x) = x$ is shown similarly.
\epf

\epoint{Remark} For a topological space $X$ with a $G$ action and fixed points $b$, $x$ let $\Path(b,x)$ denote the space of paths from $b$ to $x$. Note that $\Path(b,x)$ has a $G$ action. We can therefore define a map $\kappa$ from the fixed points of $X$ to $H^1(G, \pi_1(X,b))$ by (\ref{eqkappadef}) given above. For a $K(\pi,1)$ with $G$ action, $\kappa$ is $\pi_0$ applied to the canonical map from fixed points to homotopy fixed points.

\tpoint{Observation}\label{Gequiv_baseptchange_pi1} Let $X$ be a scheme over $k$, and let $b_1$, $b_2$ be rational base points. A choice of path $\wp$ in $\Path(b_1, b_2)$ gives an isomorphism of profinite groups $\theta: \pi_1(X_{\overline{k}}, b_2) \rightarrow \pi_1(X_{\overline{k}}, b_1)$, defined $$ \theta ( \gamma ) = \wp^{-1} \gamma \wp.$$ Note that $\theta$ is not $G_k$ equivariant. Rather, for any $g$ in $G_k$, $$g \theta (\gamma) = \kappa_{(b_1, \wp)} (b_2)^{-1} \theta(g\gamma) \kappa_{(b_1, \wp)} (b_2)$$ (cf. \ref{kappagammadefn} for the definition of $\kappa_{(b_1, \wp)} (b_2))$.

\tpoint{Observation}\label{Gkinertiapreserve} Let $X$ be a smooth curve over $k$, and let $\overline{X}$ be its smooth compactification. Suppose that $x$ is a rational point of $\overline{X} - X$. Choose a rational tangential base point $b$ at $x$. Let $\gamma$ in $\pi_1(X_{\overline{k}}, b)$ be the path determined by a small loop around the puncture at $x$. Then $\gamma$ generates the inertia group at $x$ (\cite[XIII 2.12]{sga1}), and it follows that for any $g$ in $G_k$,  $g \gamma = \gamma^{m(g)}$ for some $m(g)$ in $\Zhat$, where $\Zhat$ denotes the profinite completion of $\Z$. Furthermore, $g \gamma = \gamma^{\chi(g)}$ where $\chi: G_k \rightarrow \Zhat^*$ is the cyclotomic character. One way to see this last assertion is to note that it is sufficient to assume that $X \cup x$ is non-proper and show that the kernel of $\pi_1(X_{\overline{k}}, b)^{ab} \rightarrow \pi_1((X \cup x)_{\overline{k}}, b)^{ab}$ is $\Zhat (\chi)$. Denote this kernel by $M$. As a profinite group, $M \cong \Zhat$. $\Hom (\pi_1(X_{\overline{k}}, b)^{ab}, \Zhat)$ is the \'etale cohomology group $H^1(X_{\overline{k}}, \Zhat)$ and the analogous statement holds with $(X \cup x)_{\overline{k}}$ replacing $X_{\overline{k}}$. By the long exact sequence in cohomology of the pair $((X \cup x)_{\overline{k}}, X_{\overline{k}})$ and the purity isomorphism $$H^*((X \cup x)_{\overline{k}}, X_{\overline{k}}, \Zhat)  \cong H^{*-2}(x_{\overline{k}}, \Zhat (\chi^{-1}))= \begin{cases}\Zhat (\chi^{-1}) & \text{if $*=2$}
\\
0 &\text{otherwise,}
\end{cases}$$ we have a short exact sequence $$1 \rightarrow \Hom (\pi_1((X \cup x)_{\overline{k}}, b)^{ab}, \Zhat)\rightarrow \Hom (\pi_1(X_{\overline{k}}, b)^{ab}, \Zhat)\rightarrow \Zhat(\chi^{-1}) \rightarrow 1$$ It follows that $M \cong \Zhat(\chi)$ as $G_k$ modules as desired.

\bpoint{Galois action on $\pi_1(\proj^1_{\overline{k}} - \{\infty, a_1, a_2, \ldots, a_n\})$}\label{Gkaction_puncturedP1_bpoint} Let $b_i$ be a rational tangential base point of $\proj^1_{k} - \{\infty, a_1, a_2, \ldots, a_n \}$ at $a_i$. Let $\wp_i$ in $\Path(b_1, b_i)$ be a path from $b_1$ to $b_i$ for $i=2, 3, \ldots, n$, and let $\wp_1$ be the trivial path from $b_1$ to itself. Let $\ell_i$ in $\Path(b_i, b_i)$ be the path determined by a small loop around the puncture at $a_i$. The loops based at $b_1$ defined $$\gamma_i =\wp_i^{-1} \ell_i \wp_i $$ are free topological generators for $\pi_1(\proj^1_{\overline{k}} - \{\infty, a_1, a_2, \ldots, a_n \}, b_1)$, and it follows from \ref{Gequiv_baseptchange_pi1} and \ref{Gkinertiapreserve} that the $G_k$ action on  $$\pi_1(\proj^1_{\overline{k}} - \{\infty, a_1, a_2, \ldots, a_n\}, b_1)$$ has the form $$g \gamma_i = \mathfrak{f}_i(g)^{-1} \gamma_i^{\chi(g)} \mathfrak{f}_i(g) $$ where $\mathfrak{f}_i = \kappa_{(b_1, \wp_i)} (b_i)$ and $g$ is any element of $G_k$. 

Let $\pi$ abbreviate $\pi_1(\proj^1_{\overline{k}} - \{\infty, a_1, a_2, \ldots, a_n \}, b_1)$. Choose $v_i$ in $k^*$ for $i = 1, \ldots, n$, and suppose that $b_i$ is a rational tangential point associated to the tangent vector $a_i + v_i \epsilon$ (see \ref{tgtvecA1} for this notation).  The image of $\mathfrak{f}_i$ in $H^1(G_k, \pi^{ab})$ can be expressed in terms of the Kummer map: the basis $\{ \gamma_1, \gamma_2, \ldots, \gamma_n \}$ of $\pi^{ab}$ as a free $\Zhat(\chi)$ module determines an isomorphism $H^1(G_k, \pi^{ab}) \rightarrow H^1(G_k, \Zhat(\chi))^n$. Let $(\mathfrak{f}_i)^{ab}_j$ denote the image of $\mathfrak{f}_i$ in the $j^{th}$ factor of $H^1(G_k, \Zhat(\chi))$. Let $\kappa_j$ denote the map defined in \ref{kappadefn} for $\mathbb{A}^1 - \{a_j\}$ based at $a_j + 1$. Since the \'etale fundamental group of $\mathbb{A}^1_{\overline{k}} - \{a_j\}$ is abelian, there are canonical isomorphisms between the fundamental groups of this scheme taken with respect to different base points.  In particular, $\gamma_j$ determines an isomorphism of this fundamental group with $\Zhat(\chi)$, and $\kappa_j$ can be considered to take values in $H^1(G_k, \Zhat(\chi))$. Since the cohomology class of $\mathfrak{f}_i$ could be computed by choosing a path from $b_1$ to $b_i$ passing through  $a_j + 1$, $(\mathfrak{f}_i)^{ab}_j = \kappa_j (a_i + v_i \epsilon) - \kappa_j (a_1 + v_1 \epsilon)$. By functoriality of $\kappa$ and Lemma \ref{tgtl_pts_0_Gm}, $$\kappa_j (a + v \epsilon) =  \begin{cases}
a - a_j & \mbox{ if } a \neq a_j \\
v & \mbox{ if } a = a_j
 \end{cases}$$ Thus \begin{equation}\label{kappaabcomp} (\mathfrak{f}_i)^{ab}_j =  \begin{cases}
v_i/(a_1-a_i)  & \mbox{ if } j = i \\
(a_i - a_1)/v_1 & \mbox{ if } j = 1 \\
(a_i - a_j)/(a_1 - a_j) & \mbox{ if } j \neq 1,i \\
 \end{cases}\end{equation}

In particular, it follows that if $v_1 = a_i -a_1 = -v_i$, then the quotient of $\pi^{\Sigma}$ by $\langle \gamma_j: j \neq 1,i \rangle$ is a pro-$\Sigma$ group with $G_k$ action of the form considered in \ref{monodromy_deltan_bpoint}. (Here, $\langle \gamma_i: i \neq 1,2 \rangle$ denotes the closed normal subgroup generated by the $\gamma_i$ for $i \neq 1,2$.)  

Much more interesting information is known about the $\mathfrak{f}_i$ due to contributions of Anderson, Coleman, Deligne, Ihara, Kaneko, and Yukinari -- see for instance, \cite[6.3 Thm~p.115]{Ihara_Braids_Gal_grps}.
 
\bpoint{Restriction on $\pi_1$ sections of punctured $\proj^1$}\label{bpointrestriction_pisec} Let $X= \pmk$, and base $X$ at $0 + 1 \epsilon$, where $0 + 1 \epsilon$ denotes the tangent vector $ k[\epsilon]/\langle \epsilon^2 \rangle \rightarrow \Spec k[z] $ given by $z \mapsto \epsilon$ (cf. \ref{tgtvecA1}). Fix a positive integer $n \geq 2$ and let $\pi= \pi_1(X_{\overline{k}})^{\Sigma}$, where $\Sigma$ denotes the set of primes not dividing $n!$. By (\ref{kappaabcomp}), the presentation of $\pi_1(X_{\overline{k}})$ given in \ref{Gkaction_puncturedP1_bpoint} with its $G_k$ action is of the form (\ref{G_actionpiP1-}). Thus, the calculation of $\mu_J \delta_n$ given in \ref{monodromy_deltan_bpoint} places restrictions on the sections of $\pi_1 (X) \rightarrow G_k$. 

The $\Z^{\Sigma}(\chi)$ basis $\{ \gamma_1, \gamma_2 \}$ of $\pi^{ab}$ determines an isomorphism $$H^1(G_k, \pi^{ab}) \cong H^1(G_k, \Z^{\Sigma}(\chi))^2.$$ The quotient map $\pi/[\pi]_{n+1} \rightarrow \pi^{ab}$ therefore defines a map  $$H^1(G_k, \pi/[\pi]_{n+1}) \rightarrow H^1(G_k, \Z^{\Sigma}(\chi))^2.$$ Note that the sections of $\pi_1 (X) \rightarrow G_k$ are in natural bijection with $H^1(G_k, \pi)$. These sections map to $H^1(G_k, \Z^{\Sigma}(\chi))^2$ and the image is restricted by the following corollary of Proposition \ref{muJdeltanwithmonodromy}.

\tpoint{Corollary}\label{pmkpi1restcor}{\em The image of $H^1(G_k, \pi/[\pi]_{n+1}) \rightarrow H^1(G_k, \Z^{\Sigma}(\chi))^2$ lies in the subset of elements  $x_1 \times x_2$ such that $$\langle -x_{J(1)}, -x_{J(2)}, \ldots, -x_{J(n)} \rangle = 0$$ for all $J: \{1,2,\ldots,n \} \rightarrow \{1,2\}$ which only assume the value $2$ once.}

\bpf
An element of $H^1(G_k, \pi/[\pi]_{n+1})$ determines an element $s_n$ of $H^1(G_k, \pi/[\pi]_n)$ satisfying $\delta_n(s_n) = 0$. Applying Proposition \ref{muJdeltanwithmonodromy} to $s_n$ shows the claim.
\epf

For $n=2,3$ these restrictions are studied in \cite{PIA}. 

\epoint{Remark} Note that in the presentation of $$\pi_1(\proj^1_{\overline{k}} - \{\infty, a_1, a_2, \ldots, a_m \})$$ given in  \ref{Gkaction_puncturedP1_bpoint}, it is only possible to arrange that one of the $\mathfrak{f}_i$ for $i>1$ has image contained in the commutator subgroup, so the restrictions on $\pi_1$ sections for $\proj^1_{\overline{k}} - \{\infty, a_1, a_2, \ldots, a_m \}$ placed by Proposition \ref{muJdeltanwithmonodromy} will be pulled back from a map to $\pmk$.

\epoint{Remark}\label{pi1restremarks}  Sharifi \cite[Thm 4.3]{Sharifi} shows the vanishing of the $n^{th}$ order Massey products $\langle x,x, \ldots, x,y \rangle$ in $H^2(G_k, \Z/p^{m})$ for $x,y$ in $k^*$ such that $y$ is in the image of the norm $k(\sqrt[p^{M}]{x}) \rightarrow k$, assuming $k$ contains the $(p^M)^{th}$ roots of unity and $m \leq M-r_n$, where $r_n$ is the largest integer such that $p^{r_n}\leq n$.  Furthermore, Sharifi's methods should produce similar results under weaker hypotheses and with larger coefficient rings, although this has not been written down in detail. Sharifi's result also implies the vanishing of the Massey product $\langle y,x,x, \ldots, x, \rangle$ by formal properties of Massey products; namely, if $\langle x_1, x_2,\ldots,x_n \rangle$ is defined, then $\langle x_n, x_{n-1}, \ldots, x_1\rangle$ is defined and $$\langle x_1, x_2,\ldots,x_n\rangle = \pm \langle x_n, x_{n-1}, \ldots, x_1\rangle$$ (c.f. \cite[Thm 8]{Kraines}). This suggests redundancy among the restrictions placed by Corollary \ref{pmkpi1restcor} for $n=2$ and higher $n$.

\bpoint{Massey products in the cohomology of $G_k$} Since rational base points produce sections of $\pi_1(X) \rightarrow G_k$, applying Corollary \ref{pmkpi1restcor} to these sections produces Massey products of elements of $H^1(G_k, \Z^{\Sigma}(\chi))$ which vanish. 

We identify the elements of $H^1(G_k, \Z^{\Sigma}(\chi))^2$ corresponding to the rational base points to identify these Massey products. Let $\kappa$ denote the map of \ref{kappadefn} for $X = \pmk$ based at $0 + 1 \epsilon$ (cf. \ref{tgtvecA1}), and let $\kappa^{ab}$ denote the composition of $\kappa$ with the projection $H^1(G_k, \pi) \rightarrow H^1(G_k, \Z^{\Sigma}(\chi))^2$. For an element $x$ of $k^*$, let $x$ also denote the image in $H^1(G_k, \Z^{\Sigma}(\chi))$ of $x$ under the Kummer map.

\tpoint{Lemma}\label{kabratlbasepoints}{\em \begin{itemize} \item For $x$ in $\pmk (k) = k - \{0,1\}$, $\kappa^{ab} (x) = (x, 1-x)$. 
\item For $v$ in $k^*$, $\kappa^{ab} (1+ v \epsilon) = (1, -v)$ and $\kappa^{ab} (0 + v \epsilon) = (v ,1)$.
\item For $v$ in $k^*$, $\kappa^{ab} (\iota(0 + v \epsilon)) = (1/v ,-1/v)$, where $\iota: \pmk = \Spec k[z,\frac{1}{z}, \frac{1}{z-1}] \rightarrow \pmk$ is given by $z \mapsto \frac{1}{z}$.
\end{itemize}}

\bpf
Lemma \ref{kabratlbasepoints} follows directly from \ref{Gkaction_puncturedP1_bpoint}. More specifically, applying  (\ref{kappaabcomp}) with $a_i = x$, $a_j = 1$, $a_1 + v_1 \epsilon = 0 + 1 \epsilon$, in the cases $j = 1$ and $j \neq 1,i$ gives the formula $\kappa^{ab} (x) = (x, 1-x)$ for any $x$ in $\pmk(k)$. By (\ref{kappaabcomp}) with $a_i + v_i \epsilon= 1 + v \epsilon$, in the cases $j = 1$ and $j = i$, it follows that $\kappa^{ab} (1+ v \epsilon) = (1, -v)$ for any tangential base point $1+ v \epsilon$ at $1$. Similarly, $\kappa^{ab} (0 + v \epsilon) = (v ,1)$. Note that $\iota$ induces multiplication by $-1$ on $\pi_1(\G_{m, \overline{k}},1)$. Let $K$ denote the map of \ref{kappadefn} for $\G_m$ based at $1$. By functoriality of $\kappa$, $K(\iota(0 + v \epsilon)) = - K(0 + v \epsilon)$. By Lemma \ref{tgtl_pts_0_Gm}, $K(0 + v \epsilon) = v$. The second coordinate of $\kappa^{ab} (\iota(0 + v \epsilon))$ minus the first is $-1$ by functoriality of $\kappa$ (cf. the argument producing equation (\ref{kappaabcomp})). Thus $\kappa^{ab} (\iota(0 + v \epsilon)) = (1/v ,-1/v)$. 
\epf

\tpoint{Corollary}\label{vanishMasseyCor}{\em   Let $(X,Y)$ in $H^1(G_k, \Z^{\Sigma}(\chi))^2$ be $(x^{-1}, (1-x)^{-1})$ for $x$ in $k^* -\{1\}$, or $(x,-x)$ for $x$ in $k^*$. Then the $n^{th}$ order Massey products $$\langle X,\ldots, X, Y, X, \ldots, X \rangle$$ vanish in $H^2(G_k, \Z^{\Sigma}(\chi^2))$. Here, the Massey products have $(n-1)$ factors of $X$ and one factor of $Y$. The $Y$ can occur in any position.}

\bpf
By Lemma \ref{kabratlbasepoints}, $-(X,Y)$ is in the image of $$H^1(G_k, \pi) \rightarrow H^1(G_k, \Z^{\Sigma}(\chi))^2.$$ (Note that $-(X,Y)= (x,1-x)$ or $(x^{-1},(-x)^{-1})$.) Applying Corollary \ref{pmkpi1restcor} gives the result. 
\epf  

The vanishing of these Massey products occurs with the defining systems determined by Proposition \ref{muJdeltanwithmonodromy} and $\kappa$ applied to $x \in \pmk(k)$ or $\iota(0+ x \epsilon)$ for $x$ in $k^*$.  

\epoint{Remark} Corollary \ref{vanishMasseyCor} is also true for $(X,Y) = (x,1)$ or $(1,x)$ with $x \in k^*$ by the same proof, but this result is a formal consequence of the linearity of the Massey product \cite[Lemma 6.2.4]{Fenn}\hidden{or Vogel's thesis Prop 1.2.3}, since $1$ vanishes under the Kummer map. 

\epoint{Remark}\label{Sharifivanishing} The result of Sharifi discussed in Remark \ref{pi1restremarks} gives a different proof of the vanishing of $\langle X, X, \ldots, X, Y \rangle$ and $\langle Y, X, \ldots, X, X \rangle$ reduced mod $p^m$ when $k$ contains the $(p^M)^{th}$ roots of unity for $M \geq m + r_n$, and his methods should produce more general results as well.  They also show vanishing mod $p^m$ for more general $(X,Y)$ under his hypotheses-- see \ref{pi1restremarks}.


\bibliographystyle{DN}

\bibliography{DN}


\end{document}